\documentclass[a4paper,11pt,fullpage]{article}
\usepackage[english]{babel}

\usepackage[left=1in,right=1in,top=1in,bottom=1in]{geometry}

\usepackage{amsmath}
\usepackage{amsthm}
\usepackage{amsfonts}
\usepackage{amssymb}
\usepackage{amsxtra}
\usepackage{latexsym}
\usepackage{ifthen}
\usepackage{cite}
\usepackage{esint}

\usepackage{epic}
\usepackage{eepic}
\usepackage{xcolor}

\numberwithin{equation}{section}
\newtheorem{theorem}{Theorem}[section]
\newtheorem{lemma}{Lemma}[section]
\newtheorem{remark}{Remark}[section]
\newtheorem{definition}{Definition}[section]

\def\XXint#1#2#3{{\setbox0=\hbox{$#1{#2#3}{\int}$}
     \vcenter{\hbox{$#2#3$}}\kern-.5\wd0}}

\begin{document}

\title{Parabolic De Giorgi classes with doubly nonlinear, nonstandard growth: local boundedness under exact integrability assumptions}

\author{Simone Ciani, Eurica Henriques, Mariia O. Savchenko, Igor I. Skrypnik
}

  \maketitle

\begin{abstract} \noindent  We define a suitable class $\mathcal{PDG}$ of functions bearing unbalanced  energy estimates, that are embodied by local weak subsolutions to doubly nonlinear, double-phase, Orlicz-type and fully anisotropic operators. Yet we prove that members of $\mathcal{PDG}$ are locally bounded, under critical, sub-critical and limit growth conditions typical of singular parabolic operators, with quantitative {\it a priori} estimates that follow the lines of the pioneering work of Ladyzhenskaya, Solonnikov and Uraltseva \cite{LadSolUra}. These local bounds are new in the sub-critical cases, even for the classic $p$-Laplacean equations, since no extra-integrability condition is needed.   \end{abstract}

\pagestyle{myheadings} \thispagestyle{plain}
\markboth{  }
{On the local boundedness... }

\section{Introduction}\label{Introduction}
\noindent In \cite{DeG} Ennio De Giorgi proved that the H\"older continuity of weak solutions to 
\[ \sum_{i,j} D_i(a_{i,j}(x) D_j u ) =0,\]
for $a_{i,j}$ measurable and bounded coefficients, is a property inherent to their energy, i.e. the satisfaction of Caccioppoli inequalities, rather than to the fact that $u$ solves an equation. This concept has gained a parabolic flavour in the pioneering work \cite{LadSolUra}, where the authors identify a suitable energy class $\mathcal{B}_2(\Omega_T)$, for $\Omega_T= \Omega \times (0,T]$, consisting of functions
\[u \in C_{loc}(0,T; L^2_{loc}(\Omega)) \cap L^2_{loc}(0,T; W^{1,2}_{loc}(\Omega), \qquad \Omega \subset \subset \mathbb{R}^N,\]  
satisfying, the energy inequalities
\begin{equation}
    \begin{aligned}
        \mathcal{E}_2= \sup_{(\tau-\eta, \, \tau)} &\int_{K_{\rho}} (u-k)_{\pm}^2 \zeta^2(x,t) \, dx+ \int_{\tau-\eta}^{\tau} \int_{K_{\rho}} |D [(u-k)_{\pm} \zeta]|^2\, dxdt \\
        & \leq \int_{K_{\rho}}  (u-k)_{\pm}^2\, \zeta^2(x, \tau-\eta)\, dx + \int_{\tau-\eta}^{\tau} \int_{K_{\rho}} (u-k)_{\pm}^2\bigg[ |\zeta_t|+ |D\zeta|^2\bigg] \, dxdt  \,,       
    \end{aligned}
\end{equation} \noindent for any $k \in \mathbb{R}$ and any suitable piecewise smooth cut-off function $\zeta$ vanishing on $\partial K_{\rho}$ (see also \cite{DiB, DiBGiaVes, Giu, LadUra}  and references therein). In order to understand the extent of this definition, we observe that weak solutions to parabolic equations with measurable and bounded coefficients 
\[u_t- \text{div}(A(x,t) Du)=0,\]
with $A(x,t)= ( a_{i,j}(x,t))_{i,j=1}^N$ measurable and satisfying $\lambda |\xi|^2 \leq \xi^T  A(x,t) \xi \leq \Lambda |\xi|^2$, $0<\lambda \leq \Lambda$,  are members of $\mathcal{B}_2(\Omega_T)$. More in general the parabolic quasi-minima, as introduced by \cite{Wieser}, which are functions
\[u\in L^2_{loc}(0,T; W^{1,2}_{loc}(\Omega)), \]
satisfying for some $Q\ge 1$ the variational inequality
\begin{equation*}
    \begin{aligned}
        - \iint_{Q} u \, \phi_t \, dxdt &+ \iint_{Q} F(x,t,u, Du)\, dxdt, \qquad \quad\text{being} \qquad Q = \text{supp}(\phi)\,,\\
        &\leq Q \iint_{Q} F(x,t,u-\phi, Du-D\phi)\, dxdt \, ,
    \end{aligned}
\end{equation*} \noindent with suitable assumptions on $F$ as simply $\lambda |\xi|^2\leq F(x,t,u,\xi) \leq \Lambda |\xi|^2$, can be shown to belong to $\mathcal{B}_2(\Omega_T)$. Hence, any property derived for the members of $\mathcal{B}_2(\Omega_T)$ passes to local weak solutions, minima, and parabolic quasi-minima.\newline \noindent Now, when the growth of the operator is unbalanced, as in the case of the parabolic $p$-Laplacean
\[u_t- \text{div}(A(x,t)| Du|^{p-2} Du)=0,\]
the energy estimates $\mathcal{E}_p$ suffer from an inhomogeneity on the parabolic and elliptic energy term: 
\begin{equation} \label{p-laplacean}
    \begin{aligned}
        \mathcal{E}_p= \sup_{(\tau-\eta, \, \tau)} &\int_{K_{\rho}} (u-k)_{\pm}^2 \zeta^2(x,t) \, dx+ \int_{\tau-\eta}^{\tau} \int_{K_{\rho}} |D [(u-k)_{\pm} \zeta]|^p\, dxdt \\
        & \leq \int_{K_{\rho}}  (u-k)_{\pm}^2\, \zeta^2(x, \tau-\eta)\, dx + \int_{\tau-\eta}^{\tau} \int_{K_{\rho}} (u-k)_{\pm}^2\bigg[ |\zeta_t|+ |D\zeta|^p\bigg] \, dxdt  \,,       
    \end{aligned}
\end{equation} and still nowadays a complete regularity theory for the class of functions satisfying these inequalities is not available, since more energy inequalities tied to the structure of the equation are needed (see \cite{DiB-Annali}, \cite{DiB}) to control the oscillation, while a complete study of the local boundedness is still possible. The local boundedness of weak solutions, together with quantitative bounds on the essential suprema, is the starting point for the study of higher interior regularity such as $C^{0,\alpha}$, $C^{1,\alpha}$ interior regularity, Harnack-type inequalities and much more (see, for example \cite{LadUra, LadSolUra, DiB, DiBGiaVes, Iva}). \vskip0.1cm \noindent In this work we generalize this principle to operators that have a very wild unbalance in the energy, as for instance the doubly nonlinear parabolic equation with standard growth
\begin{equation}\label{eq4.1}
u_t- div\Big(|\nabla u^m|^{p-2} \nabla u^m\Big)=0,\quad m>0,
    \end{equation} \noindent the doubly nonlinear parabolic equation with generalized Orlicz growth $\varphi$ (see subsection \ref{ex4.2})
\begin{equation}\label{eq4.2}
u_t- div\Big(\varphi(x, t, u, |\nabla u^{m}|)\frac{\nabla u^{m}}{|\nabla u^{m}|}\Big)=0,
\end{equation} \noindent and  the doubly nonlinear anisotropic parabolic equation
\begin{equation}\label{eq4.6}
u_t-\sum\limits_{i=1}^N D_i\Big(u^{(m_i-m)(p_i-1)}|D_i u^{m}|^{p_i-2} D_i u^{m}\Big)=0\,, 
\end{equation} \noindent as well as their quasi-linear generalizations with bounded and measurable coefficients. In connection with the growth of these equations, we consider set of numbers $\{p_i\}_{i=1}^N$, $\{q_i\}_{i=1}^N$, $\{m_i\}_{i=1}^N$, $\{n_i\}_{i=1}^N$, where the $p_i$'s and $q_i$'s are linked by 
\[1<p_i\leqslant q_i\, \quad \forall i =1,2,\dots, N,\] and are meant to reflect the unbalance of the operator in its elliptic structure - such as double-phase $(p,q)$ growth, fully anisotropic, Orlicz -, wile $m_i$'s and $n_i$'s are meant to represent the unbalance offered from the point of view of the time derivative, when considering equation \eqref{eq4.1} as 
\[ (|u|^{l-1}u)_t- \Delta_p u=0\, , \qquad l=1/m\,.\] 
We define the numbers
\begin{equation}\label{exponents}
    0<m:=\min(m_1, m_2, ..., m_N, n_1, n_2, ..., n_N), \qquad q:=\max(q_1, ..., q_N),
\end{equation}\noindent the parabolic anisotropic energy spaces in a cylindrical domain $\Omega_T= \Omega \times (0,T]$ as  
\begin{equation}
\begin{aligned}
    \mathcal{V}_{loc}(\Omega_T):= \bigg{\{}  u: \Omega_T & \rightarrow \mathbb{R} \quad \text{measurable and such that for all}\, \, i=1,\dots, N,\quad \\
    &\quad u^{\frac{m_i(p_i-1)+m}{p_i}}, u^{\frac{n_i(q_i-1)+m}{p_i}}, D_i \big(u^{\frac{m_i(p_i-1)+m}{p_i}}\big)\in L^{p_i}(0, T; L^{p_i}_{loc}(\Omega))   \bigg{\}},
    \end{aligned}
\end{equation}
and the geometric configurations
\begin{equation}
    Q_{\vec{r}, \eta}(y, \tau):=K_{\vec{r}}(y)\times (\tau-\eta, \tau), \quad \text{with}  \quad K_{\vec{r}}(y):=\prod_{i=1}^N \bigg{\{} |x_i-y_i|<r_i, \bigg{\}}\,.
\end{equation}

\noindent With this notation, we can finally define the suitable parabolic De Giorgi class in the cylinder $\Omega_T$, with $\Omega \subset  \mathbb{R}^N$ open and bounded. 

\begin{definition} \label{def1}
\noindent A non-negative function $u$ belongs to De Giorgi's class $\mathcal{PDG}^+(\Omega_T, C)$ if
\[u\in C_{loc}(0, T; L^{1+m}_{loc}(\Omega))\cap \mathcal{V}_{loc}(\Omega_T),\]
 and for  every $k>0$, every cylinder $Q_{8\vec{r}, 8\eta}(y, \tau)\subset \Omega_T$, and every piecewise smooth, cutoff function $\zeta(x, t)$ vanishing on $\partial K_{\vec{r}}(y)$ and such that $0\leqslant \zeta(x, t) \leqslant 1$ the inequality 
\begin{multline}\label{eq1.1}
\mathcal{E}^m_{p}= \sup\limits_{\tau-\eta\leqslant t \leqslant \tau}\int\limits_{K_{\vec{r}}(y)}g(u^{m}, k^{m})\,\zeta^{q}\,dx+\frac{1}{C}\sum\limits_{i=1}^N  \iint\limits_{Q_{\vec{r},\eta}(y, \tau)}u^{(m_i-m)(p_i-1)}|D_i(u^{m}-k^{m})_{+}|^{p_i}\zeta^{q}\,dx\,dt\\\leqslant \int\limits_{K_{\vec{r}}(y)\times\{\tau-\eta\}}g(u^{m}, k^{m})\,\,\zeta^{q}\,dx+ C\iint\limits_{Q_{\vec{r},\eta}(y, \tau)} g(u^{m}, k^{m})\,|\zeta_t|\,dx dt+\\\quad +C\,\sum\limits_{i=1}^N \iint\limits_{Q_{\vec{r},\eta}(y, \tau)}u^{(m_i-m)(p_i-1)}(u^{m}-k^{m})_{+}^{p_i}|D_i \zeta|^{p_i}\,dx\,dt+ \\+C\,\sum\limits_{i=1}^N \iint\limits_{Q_{\vec{r},\eta}(y, \tau)}u^{(n_i-m)(q_i-1)}(u^{m}-k^{m})_{+}^{q_i}|D_i \zeta|^{q_i}\,dx\,dt,
\end{multline}
is valid with a constant $C>0$ depending only on the data\footnote{We say that a constant $C$ depends only on the data if it depends only on $\{p_i, q_i, n_i, m_i, N\}$}, and 
\[g(u^{m}, k^{m}):=\frac{1}{m}\int\limits^{u^{m}}_{k^{m}} z^{\frac{1}{m}-1} (z-k^{m})_+\,dz.\]
\end{definition}
\noindent 

\begin{remark} In the case of the standard $p$-Laplacean equation $u_t-\Delta_p u=0$, i.e. in \eqref{eq4.1} with $m=1$, the function $g(u,k)$ above reduces to the classic parabolic term  $g(u,k)= (u-k)_+^2 /2$.
\end{remark}

\begin{remark}
    Since $u \in \mathcal{V}_{loc}(\Omega_T)$, then all the integral quantities in \eqref{eq1.1} are convergent, and the chain rule
    \[\bigg|D_i \bigg( [u^m]^{(\frac{m_i(p_i-1)+m}{p_i m})}\bigg)\bigg|^{p_i}= \gamma_i  |D_i [u^m]|^{p_i} \, u^{(m_i-m)(p_i-1)}\,,\qquad \quad \gamma_i=\bigg(\frac{m_i(p_i-1)+m}{p_i m}  \bigg)\,,\] will be tacitly assumed in the membership $u \in \mathcal{V}_{loc}(\Omega_T)$.
\end{remark}

\subsection*{Novelty and Significance} Motivation of this paper relies on the relatively recent research in \cite{Ok}, where the regularity results for solutions of double phase elliptic equation with $(p,q)$ growth were proved under additional assumptions that 
\begin{equation}\label{eq1.2}
q>\frac{Np}{N-p},\quad p<N\quad \text{and}\quad u\in L^s_{loc}(\Omega)\quad \text{with}\quad s>\frac{N}{p}(q-p)>\frac{N p}{N-p}.
\end{equation}
\noindent The boundedness of solutions to elliptic  equations with $(p,q)$ growth in the case $q\leqslant \frac{N p}{N-p}$ (as well, to parabolic equations under condition $\frac{2 N}{N+2}<p\leqslant q \leqslant p(1+\frac{2}{N})$) starting with the papers of Kolodii \cite{Kol1, Kol2} is well known, see, for example \cite{BiaCupMas, BocMarSbo, CarGaoGioLeo, Cia, CupMarMas, CupLeoMas, FusSbo1, FusSbo2, Str, MinXit, Sin}.
At the same time,  examples constructed by Marcellini \cite{Mar} and Giaquinta \cite{Gia} show that there exist unbounded solutions, provided that $q>\frac{N p}{N-p}$. In this case, additional assumptions are needed on the solution for its local boundedness, that are extra-integral conditions as \eqref{eq1.2}. \newline \noindent Similar questions arise for parabolic equations with nonstandard growth, but, unlike elliptic equations, in this case two special exponents show up: one of them (see case $\Lambda>1$ n the next Section) is responsible for the behavior of subsolutions in the degenerate case, and the second (see case $\Lambda<1$ n the next Section) in the singular case. Surprisingly, even in the case of standard growth, i.e. referring to $\mathcal{PDG}^+(\Omega_T)$ with $0<m = m_i = n_i$,
$p = p_i = q_i$, $i = 1, ..., N$ and $(m+1)p+N(m(p-1)-1)=0$ (if $m=1$ that is $p=\frac{2N}{N+2}$, or if $p=2$ that is $m=\frac{N-2}{N+2}$, $N\geqslant 3$), local boundedness was known only under the additional condition that $u\in L^s_{loc}(\Omega_T)$ with sufficiently large $s$ (see, for example \cite{DiB, BogDuzGiaLiaSch, BonVaz, BonIagVaz, DiBGiaVes}, \cite{Henriques1}, \cite{Henriques2}). As Theorem \ref{th1.1} shows, this condition can be removed. \newline \noindent About the known counter-examples (see for instance \cite{BogDuzGiaLiaSch}), simple computations show that the known unbounded solutions in this borderline case do not belong to the corresponding Sobolev spaces defining the solutions, and so, the corresponding unbounded solutions fail to be weak solutions. Here we note also that in the sub-subcritical case the quantitative point-wise estimate is usually proved under the additional condition that  $u$ is qualitatively locally bounded (see \cite{DiB, BogDuzGiaLiaSch, BonVaz, BonIagVaz, DiBGiaVes} for the standard case). Theorem \ref{th1.2} shows that this condition can be removed too.

\subsection*{Main Results} In order to state our main results we need to define numbers $\{\lambda_i\}_{i=1}^N$, $\{\Lambda_i\}_{i=1}^N$ by
\[\lambda_i:=m_i(p_i-1)\leqslant \Lambda_i:=n_i(q_i-1)\qquad \forall i=1, ..., N,\]
relative to the set of numbers $\{p_i\}_{i=1}^N$, $\{q_i\}_{i=1}^N$, $\{m_i\}_{i=1}^N$, $\{n_i\}_{i=1}^N$ of $\mathcal{PDG}^+(\Omega_T)$. Then we define the special indexes 
\[
\Lambda:=\max(\Lambda_1, ..., \Lambda_N),\quad
\frac{1}{p}:=\frac{1}{N}\sum\limits_{i=1}^N \frac{1}{p_i},\quad  \Big|\frac{\lambda}{p}\Big|:=\frac{1}{N}\sum\limits_{i=1}^N \frac{\lambda_i}{p_i},\quad p<N.
\]

\noindent When $\Lambda <1$ we say that the equation has {\it singular} behavior, while when $\Lambda >1$ we say that the equation has {\it degenerate} behavior. Clearly, we refer to the case $\Lambda=1$ as the {\it limiting} case.

\begin{theorem}\label{th1.1}
Let $u \in \mathcal{PDG}^+(\Omega_T)$ and let
\begin{equation}\label{eq1.3}
\mathcal{L}:=\max (1, \Lambda)\leqslant p\Big|\frac{\lambda}{p}\Big|+(m+1)\frac{p}{N}=:\mathcal{M},
\end{equation}
then $u$ is locally bounded. In particular, if
\begin{equation*}
\mathcal{L} < p\Big|\frac{\lambda}{p}\Big|+(m+1)\frac{p}{N},
\end{equation*}
then for  any cylinder $Q_{8 \vec{\rho}, 8\theta}(x_0, t_0)\subset \Omega_T$ the following estimate holds true
\begin{equation}\label{general-formula}
    \sup\limits_{Q_{\frac{\vec{\rho}}{2}, \frac{\theta}{2}}(x_0, t_0)} u\leqslant \gamma \bigg(\mathcal{H}(\theta, \vec{\rho} )\iint\limits_{Q_{\vec{\rho}, \theta}(x_0, t_0)} u^{m+\mathcal{L}}\,dx dt\bigg)^{\frac{p}{(\mathcal{M}-\mathcal{L})N}} + \mathcal{R}(\theta, \vec{\rho})\,,
\end{equation} for a constant $\gamma>0$ that depends only on the data, while the functions $\mathcal{H}(\theta, \vec{\rho} ), \mathcal{R}(\theta, \vec{\rho})$ are defined, respectively, in \eqref{eq3.4} and \eqref{eq3.5}. 

\end{theorem}
\noindent Now, in order to formulate the respective of Theorem \ref{th1.1} in the sub-critical case $\mathcal{L} \ge\mathcal{M}$, we need to define the non-degeneracy number
\begin{equation}\label{eq1.7}
\varkappa_s:=p(s+1-\mathcal{L})+N(p\Big|\frac{\lambda}{p}\Big|-\mathcal{L})=
\begin{cases}
p s+N(p\big|\frac{\lambda}{p}\big|-1),\quad \text{if}\quad \Lambda< 1,\\
p(s+1-\Lambda)+N(p\big|\frac{\lambda}{p}\big|-\Lambda),\quad \text{if}\quad \Lambda> 1.
\end{cases}
\end{equation}

\begin{theorem}\label{th1.2}
Let $u\in \mathcal{PDG}^+(\Omega_T)$, let $\mathcal{L} \ge \mathcal{M}$
% i.e. 
% \begin{equation}\label{eq1.8}
% \mathcal{L} > \mathcal{M},
% \end{equation}
\noindent and let $s> m+1$ be any number such that
\begin{equation}\label{eq1.9}
\varkappa_s >0,\quad \text{i.e.} \quad s> \mathcal{L}-1 +\frac{N}{p}(\mathcal{L}- p\Big|\frac{\lambda}{p}\Big|).
\end{equation}
If we assume additionally  that $u\in L^s_{loc}(\Omega_T)$,  then 
$u$ is locally bounded, and moreover, the following inequality holds
\begin{equation}\label{general-formula-2}
    \sup\limits_{Q_{\frac{\vec{\rho}}{2}, \frac{\theta}{2}}(x_0, t_0)} u\leqslant \gamma \bigg(\mathcal{H}(\theta, \vec{\rho})\iint\limits_{Q_{\vec{\rho}, \theta}(x_0, t_0)} u^{s}\,dx dt\bigg)^{\frac{p}{\varkappa_s}} +\mathcal{R}(\theta, \vec{\rho})\,,
\end{equation}
% \begin{multline}\label{eq1.10}
% \sup\limits_{Q_{\frac{\vec{\rho}}{2}, \frac{\theta}{2}}(x_0, t_0)} u\leqslant \gamma \bigg(\frac{1}{\theta^{\frac{N+p}{p}}}\iint\limits_{Q_{\vec{\rho}, \theta}(x_0, t_0)} u^{s}\,dx dt\bigg)^{\frac{p}{\varkappa_s}} +\\+
% \gamma \bigg[\sum\limits_{i=1}^N\Big(\frac{\theta}{\rho^{p_i}_i}\Big)^{\frac{1}{1-\lambda_i}}+ \sum\limits_{i=1}^N\Big(\frac{\theta }{\rho^{q_i}}\Big)^{\frac{1}{1-\Lambda_i}}\bigg],\quad \text{if}\quad \Lambda<1=\mathcal{L},
% \end{multline}
% and
% \begin{multline}\label{eq1.11}
% \sup\limits_{Q_{\frac{\vec{\rho}}{2}, \frac{\theta}{2}}(x_0, t_0)} u\leqslant \gamma \bigg(\Big[\sum\limits_{i:\lambda_i=\Lambda} \frac{1}{\rho^{p_i}_i}+\sum\limits_{i:\Lambda_i=\Lambda} \frac{1}{\rho^{q_i}_i}\Big]^{\frac{N+p}{p}}\iint\limits_{Q_{\vec{\rho}, \theta}(x_0, t_0)} u^{s}\,dx dt\bigg)^{\frac{p}{\varkappa_s}} +\\+\gamma \bigg[\sum\limits_{i:\lambda_i=\Lambda}\Big(\frac{\rho^{p_i}_i}{\theta}\Big)^{\frac{1}{\Lambda-1}}+
% \sum\limits_{i:\Lambda_i=\Lambda}\Big(\frac{\rho^{q_i}_i}{\theta}\Big)^{\frac{1}{\Lambda-1}}
% +\sum\limits_{i:\lambda_i<\Lambda}\Big(\frac{1}{\rho^{p_i}_i}\sum\limits_{l:\lambda_l=\Lambda}\rho^{p_l}_l\Big)
% ^{\frac{1}{\Lambda-\lambda_i}}+\\+\sum\limits_{i:\Lambda_i<\Lambda}\Big(\frac{1}{\rho^{q_i}_i}
% \sum\limits_{l:\Lambda_l=\Lambda}\rho^{q_l}_l\Big)
% ^{\frac{1}{\Lambda-\Lambda_i}}\bigg],\quad \text{if} \quad 1< \Lambda=\mathcal{L},
% \end{multline}
for a  constant $\gamma>0$ depending only on the data.
% Likewise, in the limiting case $\mathcal{L} =\mathcal{M}$, when the additional information $u\in L^s_{loc}(\Omega_T)$ is available with $s$ satisfying \eqref{eq1.9}, then \eqref{general-formula-2} holds true.
\end{theorem}

\begin{remark}
The condition $\varkappa_s>0$ allows us to control the behavior of $u$ both in the singular case ($\Lambda<1$) and in the degenerate case ($\Lambda>1$). Our  approach allows us to prove the $L^\infty_{loc}(\Omega_T)$ bound in both of these cases simultaneously.\end{remark} \noindent By simple computations on the exponents one sees that for the $p$-Laplacean equations \eqref{p-laplacean} the condition $\mathcal{L}<\mathcal{M}$ corresponds to $p>(N+2)/(2N)$, so that, coherently with the study of \cite{DiB}, we refer to Theorem \eqref{th1.1} as the {\it super-critical} case, and Theorem \ref{th1.2} as the {\it critical/sub-critical case}.

\subsubsection*{Structure of the Proof} The proofs of Theorems \ref{th1.1}-\ref{th1.2} are based on the method of De Giorgi \cite{DeG}, for an exhaustive overview on the subject, see, for example \cite{Giu, LadUra, LadSolUra, DiB, DiBGiaVes}. Although the strategy for establishing boundedness goes as in the standard $p$-growth, we had to overcome some
difficulties for the presence of $(m_i, p_i)$- $(n_i, q_i)$-growth that cause the unbalanced   energy inequality \eqref{eq1.1} defining $\mathcal{PDG}^+(\Omega_T)$. In the limiting case $\mathcal{L} = \mathcal{M}$ we also explore the ideas from \cite{BocMarSbo, CupMarMas, FusSbo1, FusSbo2}.\vskip0.2cm\noindent 
\subsubsection*{Comparison with known literature} Theorem \ref{th1.1} generalizes Ok's result \cite{Ok} for  elliptic equations with $(p, q)$ growth and, in addition, refines known results in the so-called  sub-subcritical case for doubly nonlinear parabolic equations with standard $(m, p)$-growth, \cite{BogDuzGiaLiaSch}. In particular, if 
\[(m+1) p+N (m(p-1)-1)=0,\quad \qquad  m>0,\]
we obtain the local boundedness of non-negative subsolutions without the assumption that \\$u\in L^s_{loc}(\Omega_T)$ with some sufficiently large $s>1$.  Our results also cover the new cases of parabolic equations with nonstandard  $(m,p)-(n,q)$ and anisotropic $(m_i, p_i)$  growth, generalizing the parabolic studies of \cite{CH}, \cite{CHS},\cite{Porzio2}, \cite{DMV}, and literature therein. See subsection \ref{ex4.3} for more details.
\vskip0.2cm \noindent {\bf Structure of the paper:} In Section $2$  we collect some preliminary and technical properties,  Section $3$  is devoted to the proofs theorems \ref{th1.1}-\ref{th1.2}. Finally, Section $4$ contains several examples illustrating that solutions to equations \eqref{eq4.1}-\eqref{eq4.2}-\eqref{eq4.6} belong to $\mathcal{PDG}^+(\Omega_T)$, and therefore Theorems \ref{th1.1}, \ref{th1.2} apply to them.

\section{Auxiliary material}\label{Sec2}
\noindent In this brief section we collect our main tools of the trade. First, we start with a Lemma that allows to estimate more precisely function $g$ in \eqref{eq1.1} of Definition \ref{def1}.
For $u$, $k$, $m\in \mathbb{R}_+$ set
$$g_{\pm}(u, k):=\pm\frac{1}{m}\int\limits^u_k z^{\frac{1}{m}-1} (z-k)_{\pm}\,dz.$$
Then, truncations $g_{\pm}(u, k)$ are bounded above and below.
\begin{lemma}[e.g. \cite{BogDuzGiaLiaSch}, Lemma $3.2$]\label{lem2.1}
There exists constant $\gamma >0$ depending only on $m$ such that
\begin{equation*}
\frac{1}{\gamma}\,(k+u)^{\frac{1}{m}-1}(u-k)^2_{\pm}\leqslant g_{\pm}(u, k)\leqslant \gamma\, (k+u)^{\frac{1}{m}-1}(u-k)^2_{\pm}.
\end{equation*}
\end{lemma} \noindent The next Lemma is called in the literature The Geometric Convergence Lemma, it is a classic tool of for {\it a priori} estimates for which we refer to \cite{DiB}, Chapter $2$.
\begin{lemma}\label{lem2.2}
Let $\{y_j\}_{j\in \mathbb{N}}$ be a sequence of positive numbers satisfying the recursive inequalities
$$y_{j+1}\leqslant A\,B^j\,y^{1+\delta}_j,\quad j=0, 1, ...,$$
where $A$, $B>1$ and $\delta>0$ are given numbers.  If 
$$y_0\leqslant A^{-\frac{1}{\delta}}\,B^{-\frac{1}{\delta^2}},$$
then $\lim\limits_{j\rightarrow \infty} y_j=0.$
\end{lemma}
\noindent Finally, we will need the following anisotropic embedding at its full strength.
\begin{lemma}[e.g. \cite{Tro} or \cite{DMV} Prop. 2.3.] \label{lem2.3}
Let $\Omega$ be an open
bounded set in $\mathbb{R}^N$, $u\in W^{1,1}_0(\Omega)$ and let also
$$\sum\limits_{i=1}^N \int\limits_{\Omega} |D_i u^{\alpha_i}|^{p_i} dx <\infty,\quad \alpha_i>0,\quad p<N.$$ 
Then there exists $c$, depending on $N, p_1, ..., p_N, \alpha_1, ..., \alpha_N$ such that
\begin{equation*}
\Big(\int\limits_{\Omega}|u|^{\bar{p}} dx\Big)^{\frac{N-p}{N}}\leqslant c \prod\limits_{i=1}^N \Big(\int\limits_{\Omega}|D_i u^{\alpha_i}|^{p_i} dx \Big)^{\frac{p}{N p_i}},\quad \bar{p}=\frac{N \alpha p}{N-p},\quad \alpha=\frac{1}{N}\sum\limits_{i=1}^N \alpha_i.
\end{equation*}
\end{lemma}

\subsection{Notation} We refer to the parameters $C, N, m_1, ..., m_N, n_1, ..., n_N, p_1, ..., p_N, q_1, ..., q_N$ as our structural data, and
we write $\gamma$ if it can be quantitatively determined a priori only in terms of the above quantities. The generic constant $\gamma$ may change from line to line.

\section{Local Boundedness: Proof of Theorems \ref{th1.1}, \ref{th1.2}}\label{Sec.3}
\noindent The first step is a common energy estimate. We set up the iterative geometry: let us fix $(x_0, t_0)\in \Omega_T$ such that $Q_{8\vec{\rho}, 8\theta}(x_0, t_0)\subset \Omega_T$. For $k>0$ to be determined, we define the increasing sequence of levels  
\[k_j^{m}=k^{m}-\frac{k^{m}}{2^{j+1}}, \qquad \bar{k}^{m}_j:=\frac{1}{2}(k^{m}_j+k^{m}_{j+1}), \qquad [k^{'}_j]^{m}:=\frac{1}{2}(\bar{k}^{m}_j+k^{m}_{j+1}),\]
and the decreasing geometric sequences
\[\vec{r}_j:= \frac{\vec{\rho}}{2}\bigg(1+\frac{1}{2^{j+1}}\bigg), \qquad \eta_j:=\frac{\theta}{2}\bigg(1 +\frac{1}{2^{j+1}}\bigg),\]
\[K_j:=K_{\vec{r}_j}(x_0), \qquad Q_j:=K_j\times I_j, \qquad I_j:= (t_0-\eta_j, t_0)\,.\]
\vskip0.2cm \noindent Now, let $\zeta^{(1)}_j(x)\in C^1_0(K_j)$ be a cut-off function between $K_j$ and $K_{j+1}$, i.e.
\[\zeta^{(1)}_j(x)=1 \quad \text{in} \quad K_{j+1}, \quad 0\leqslant \zeta^{(1)}_j(x)\leqslant 1,\] and therefore obliged to satisfy the decays 
\[|D_i \zeta^{(1)}_{j}(x)|\leqslant \gamma \frac{2^j}{\rho_i}, \qquad \forall i=1, ..., N\,. \] Moreover, in order to localize with respect to time, let us consider $\zeta^{(2)}_j(t)\in C^1(\mathbb{R})$ such that 
\[\zeta^{(2)}_j(t)=0 \quad \text{for} \quad t\leqslant t_0-\eta_j, \quad\zeta^{(2)}_j(t)=1 \quad \text{for}\quad t\geqslant t_0-\eta_{j+1},\]
\[1\leqslant \zeta^{(2)}_j(t)\leqslant 1, \quad |\zeta^{(2)}_{j, t}(t)|\leqslant \gamma \frac{2^j}{\theta},\]
and finally set $\zeta_j:=\zeta^{(1)}_j \zeta^{(2)}_j$.
\noindent In order to control the parabolic energy terms, we observe that on the set $\{u\geqslant k_{j+1}\}$ one has
\[1\leqslant \frac{u^{m}+[k^{'}_j]^{m}}{u^{m}- [k^{'}_j]^{m}}\leqslant \frac{2u^{m}}{u^{m}- [k^{'}_j]^{m}}\leqslant \frac{2k^{m}_{j+1}}{k^{m}_{j+1}-[k^{'}_j]^{m}}\leqslant \gamma 2^{j\gamma},\]
and hence Lemma \ref{lem2.1} yields
\begin{equation}\label{eq3.1}
g(u^{m}, [k^{'}_j]^{m})\geqslant \frac{1}{\gamma 2^{j\gamma}}(u^{m}-k^{m}_{j+1})^{1+\frac{1}{m}},\quad \text{if}\quad u\geqslant k_{j+1}.
\end{equation} 
Moreover, on the set $\{u\geqslant k^{'}_j\}$
$$1\leqslant \frac{u^{m}+[k^{'}_j]^{m}}{u^{m}- \bar{k}^{m}_j}\leqslant \frac{2u^{m}}{u^{m}- \bar{k}^{m}_j}\leqslant \frac{2 [k^{'}_j]^{m}}{[k^{'}_j]^{m} -\bar{k}^{m}_j}\leqslant \gamma 2^{j\gamma},$$
and hence, using the fact that $(u^{m}-[k^{'}_j]^{m})_+\leqslant (u^{m}-\bar{k}^{m}_j)_+$, from Lemma \ref{lem2.1} we get
\begin{equation}\label{eq3.2}
g(u^{m}, [k^{'}_j]^{m})\leqslant \gamma 2^{j\gamma}(u^{m}-\bar{k}^{m}_j)^{1+\frac{1}{m}}_+.
\end{equation}
Therefore, setting
$$\alpha_i:=\frac{\lambda_i+m}{m\,p_i},\quad i=1, ..., N,$$
by \eqref{eq3.1}, \eqref{eq3.2},  inequality \eqref{eq1.1} can be rewritten as
\begin{multline}\label{eq3.3}
\sup\limits_{t\in I_j}\int\limits_{K_j}(u^{m}-k^{m}_{j+1})^{1+\frac{1}{m}}_{+}\zeta^{q}_j\,dx+
\sum\limits_{i=1}^N\iint\limits_{Q_{j}}|D_i (u^{m}- k^{m}_{j+1})^{\alpha_i}_{+}|^{p_i}\zeta^{q}_j\,dx\,dt\\\leqslant 
\frac{\gamma 2^{j\gamma}}{\theta}\iint\limits_{Q_{j}}(u^{m}-\bar{k}^{m}_{j})^{1+\frac{1}{m}}_{+}\,\,dx dt+\\+
\sum\limits_{i=1}^N \frac{\gamma 2^{j\gamma}}{\rho^{p_i}_i}\,\iint\limits_{Q_{j}}u^{(m_i-m)(p_i-1)}(u^{m}-\bar{k}^{m}_{j})_{+}^{p_i}\,dx\,dt+\\+\sum\limits_{i=1}^N\frac{\gamma 2^{j\gamma}}{\rho^{q_i}_i}\,\iint\limits_{Q_{j}}u^{(n_i-m)(q_i-1)}(u^{m}-\bar{k}^{m}_{j})_{+}^{q_i}\,dx\,dt
\\\leqslant
\gamma 2^{j\gamma}\bigg[\frac{k^{1-\mathcal{L}}}{\theta}+\sum\limits_{i=1}^N\frac{k^{\lambda_i-\mathcal{L}}}{\rho^{p_i}_i}+\sum\limits_{i=1}^N
\frac{k^{\Lambda_i-\mathcal{L}}}{\rho^{q_i}_i} \bigg]\,\iint\limits_{Q_{j}}(u^{m}- k^{m}_{j})_{+}^{1+\frac{\mathcal{L}}{m}}\,dx\,dt.
\end{multline}
By the simple arithmetic fact that $\sum\limits_{i=1}^{S}a_i\geqslant \Big(\sum\limits_{i=1}^S\frac{1}{a_i}\Big)^{-1}$ for any $S\geqslant 1$ and $a_i>0$, the  terms in brackets in \eqref{eq3.3} are estimated by stipulating to take
\begin{equation}\label{eq3.4}
k\geqslant \mathcal{R}(\theta, \vec{\rho}):=
\begin{cases}
\sum\limits_{i=1}^N\big(\frac{\theta}{\rho^{p_i}_i}\big)^{\frac{1}{1-\lambda_i}}+ \sum\limits_{i=1}^N \big(\frac{\theta}{\rho^{q_i}_i}\big)^{\frac{1}{1-\Lambda_i}},\quad \text{if}\quad \Lambda<1=\mathcal{L},\\
\sum\limits_{i:\lambda_i<1}\big(\frac{\theta}{\rho^{p_i}_i}\big)^{\frac{1}{1-\lambda_i}}+
\sum\limits_{i:\Lambda_i<1}\big(\frac{\theta}{\rho^{q_i}_i}\big)^{\frac{1}{1-\Lambda_i}},\quad \text{if}\quad \Lambda=1=\mathcal{L},\\
\sum\limits_{i:\lambda_i=\Lambda}\big(\frac{\rho^{p_i}_i}{\theta}\big)^{\frac{1}{\Lambda-1}}+
\sum\limits_{i:\Lambda_i=\Lambda}\big(\frac{\rho^{q_i}_i}{\theta}\big)^{\frac{1}{\Lambda-1}}
+\sum\limits_{i:\lambda_i<\Lambda}\Big(\frac{1}{\rho^{p_i}_i}\sum\limits_{l:\lambda_l=\Lambda}\rho^{p_l}_l\Big)
^{\frac{1}{\Lambda-\lambda_i}}+\\+\sum\limits_{i:\Lambda_i<\Lambda}\Big(\frac{1}{\rho^{q_i}_i}
\sum\limits_{l:\Lambda_l=\Lambda}\rho^{q_l}_l\Big)
^{\frac{1}{\Lambda-\Lambda_i}},\quad \text{if} \quad 1< \Lambda=\mathcal{L}.
\end{cases}
\end{equation}
Therefore, setting
\begin{equation}\label{eq3.5}
\mathcal{H}(\theta, \vec{\rho}):=
\begin{cases}
\frac{1}{\theta},\quad \text{if}\quad \Lambda<1,\\
\frac{1}{\theta}+\sum\limits_{i:\lambda_i=1} \frac{1}{\rho^{p_i}_i}+\sum\limits_{i:\Lambda_i=1} \frac{1}{\rho^{q_i}_i},\quad \text{if}\quad \Lambda=1,\\
\sum\limits_{i:\lambda_i=\Lambda} \frac{1}{\rho^{p_i}_i}+\sum\limits_{i:\Lambda_i=\Lambda} \frac{1}{\rho^{q_i}_i},\quad \text{if}\quad \Lambda>1,
\end{cases}
\end{equation}
the energy inequality \eqref{eq3.3} gets simplified into
\begin{multline}\label{eq3.6}
\sup\limits_{t\in I_j}\int\limits_{B_j}(u^{m}-k^{m}_{j+1})^{1+\frac{1}{m}}_{+}\zeta^{q}_j\,dx+
\sum\limits_{i=1}^N\iint\limits_{Q_{j}}|D_i (u^{m}- k^{m}_{j+1})^{\alpha_i}_{+}|^{p_i}\zeta^{q}_j\,dx\,dt\\\leqslant 
\gamma 2^{j\gamma}\,\mathcal{H}(\theta, \vec{\rho})\,\iint\limits_{Q_{j}}(u^{m}- k^{m}_{j})_{+}^{1+\frac{\mathcal{L}}{m}}\,dx\,dt.
\end{multline}
\noindent Here the proof of the two Theorems splits.
\subsection{Proof of Theorem \ref{th1.1}}

\subsubsection{Case  $\mathcal{L}< \mathcal{M}$}
Choose  
\begin{equation*}
\begin{cases}
\alpha:=\frac{1}{N}\sum\limits_{i=1}^N\alpha_i=\frac{1}{p}\big(1+\frac{p}{m}\Big|\frac{\lambda}{p}\Big|\big),\\
 p_*:= \alpha p+(1+\frac{1}{m})\frac{p}{N}> 1+\frac{\mathcal{L}}{m}.
 \end{cases}
\end{equation*}
Let us set
\[\alpha^-:=\min\limits_{1\leqslant i\leqslant N}\alpha_i, \quad p^-:=\min\limits_{1\leqslant i\leqslant N} p_i, \quad \gamma_0:=q\bigg(1+\frac{1}{p^- \alpha^-}\bigg)\,.\] \vskip0.1cm \noindent By H\"{o}lder's inequality and  Lemma \ref{lem2.3}, from \eqref{eq3.6} we obtain
\begin{equation}\label{eq3.7}
\begin{aligned}
\iint\limits_{Q_{j+1}}&(u^{m}-k^{m}_{j+1})^{1+\frac{\mathcal{L}}{m}}_{+}\,\,dx dt\\
&\leqslant \Big(\iint\limits_{Q_{j}}\big[(u^{m}-k^{m}_{j+1})_{+}\zeta^{\gamma_0}_j\big]^{p_*}\,\,dx dt\Big)^{\frac{1+\frac{\mathcal{L}}{m}}{p_*}} |Q_j\cap\{u\geqslant k_{j+1}\}|^{1-\frac{1+\frac{\mathcal{L}}{m}}{p_*}}\\
&\leqslant \gamma \Big(\sup\limits_{t\in I_j}\int\limits_{B_j}(u^{m}-k^{m}_j)^{1+\frac{1}{m}}_+\zeta^{q}_j\,dx\Big)^{\frac{p(1+\frac{\mathcal{L}}{m})}{N p_*}}  \times \\
& \qquad \times \Big(\int\limits_{I_j}\Big[\int\limits_{B_{j}}\big[(u^{m}-k^{m}_{j+1})_{+}\zeta^{\frac{q}{p^- \alpha^-}}_j\big]^{\frac{N \alpha p}{N-p}}\,\,dx \Big]^{\frac{N-p}{N}}\,dt\Big)^{\frac{1+\frac{\mathcal{L}}{m}}{p_*}}|Q_j\cap\{u\geqslant k_{j+1}\}|^{1-\frac{1+\frac{\mathcal{L}}{m}}{p_*}}\\
&\leqslant \gamma \Big(\sup\limits_{t\in I_j}\int\limits_{B_j}(u^{m}-k^{m}_j)^{1+\frac{1}{m}}_+\zeta^{q}_j\,dx\Big)^{\frac{p(1+\frac{\mathcal{L}}{m})}{N p_*}}\times \\
& \qquad \times \Big(\sum\limits_{i=1}^N\iint\limits_{Q_{j}}\big|D_i \big[(u^{m}-k^{m}_{j+1})^{\alpha_i}_{+}\zeta^{\frac{q}{p^-}}_j\big]\big|^{p_i}\,\,dx dt \Big)^{\frac{1+\frac{\mathcal{L}}{m}}{p_*}}|Q_j\cap\{u\geqslant k_{j+1}\}|^{1-\frac{1+\frac{\mathcal{L}}{m}}{p_*}},
\end{aligned}
\end{equation} \noindent and by the common energy estimates \eqref{eq3.6} we obtain 
\begin{equation}
    \begin{aligned}
        \iint\limits_{Q_{j+1}}(u^{m}&-k^{m}_{j+1})^{1+\frac{\mathcal{L}}{m}}_{+}\,\, dx dt\\
        & \leqslant \gamma 2^{j\gamma}\,k^{-(m+\mathcal{L}) (1-\frac{1+\frac{\mathcal{L}}{m}}{p_*})}\Big(\iint\limits_{Q_{j}}(u^{m}-k^{m}_{j})^{1+\frac{\mathcal{L}}{m}}_{+}\,\,dx dt\Big)^{1-\frac{1+\frac{\mathcal{L}}{m}}{p_*}}
         \times\\
         & \qquad  \quad \times \Big[\mathcal{H}(\theta, \vec{\rho})\iint\limits_{Q_{j}}(u^{m}-k^{m}_{j})^{1+\frac{\mathcal{L}}{m}}_{+}\,\,dx dt\Big]^{\frac{(1+\frac{\mathcal{L}}{m}) (N+p)}{p_* N}}.
    \end{aligned}
\end{equation}
Therefore, if we define 
$$y_j:=\iint\limits_{Q_j}(u^{m}-k^{m}_{j})^{1+\frac{\mathcal{L}}{m}}_{+}\,\,dx dt,$$
then from \eqref{eq3.7} we obtain 
\begin{equation}\label{eq3.8}
y_{j+1}\leqslant \gamma 2^{j\gamma}k^{-(m+\mathcal{L})(1-\frac{1+\frac{\mathcal{L}}{m}}{p_*})}\big[\mathcal{H}(\theta, \vec{\rho})\big]^{\frac{(1+\frac{\mathcal{L}}{m})(N+p)}{p_* N}}\, y^{1+\frac{(1+\frac{\mathcal{L}}{m})p}{p_* N}}_j.
\end{equation}
By Lemma \ref{lem2.2} $\lim\limits_{j\rightarrow \infty} y_j=0$, provided that $k$ is chosen to satisfy
\begin{equation}\label{eq3.9}
k^{\frac{(m+1)p+N(p|\frac{\lambda}{p}|-\mathcal{L})}{p}}\geqslant \gamma \big[\mathcal{H}(\theta, \vec{\rho})\big]^{\frac{N+p}{p}}\, y_0.
\end{equation}
From this, taking into account \eqref{eq3.4} and \eqref{eq3.5}, we arrive at the required estimate \eqref{general-formula}.

\begin{remark}[Explicit formulas] In the singular case  $\Lambda<1=\mathcal{L}$, formula \eqref{general-formula} takes the explicit shape 
\begin{multline}\label{eq1.4}
\sup\limits_{Q_{\frac{\vec{\rho}}{2}, \frac{\theta}{2}}(x_0, t_0)} u\leqslant \gamma \bigg(\frac{1}{\theta^{\frac{N+p}{p}}}\iint\limits_{Q_{\vec{\rho}, \theta}(x_0, t_0)} u^{m+1}\,dx dt\bigg)^{\frac{p}{(m+1)p+N(p|\frac{\lambda}{p}|-1)}} +\\+\gamma \bigg[\sum\limits_{i=1}^N\Big(\frac{\theta}{\rho^{p_i}_i}\Big)^{\frac{1}{1-\lambda_i}}+ \sum\limits_{i=1}^N\Big(\frac{\theta }{\rho^{q_i}}\Big)^{\frac{1}{1-\Lambda_i}}\bigg], 
\end{multline} \noindent while in the limiting case $\Lambda=1=\mathcal{L}$ the expression \eqref{general-formula} is 
\begin{multline}\label{eq1.5}
\sup\limits_{Q_{\frac{\vec{\rho}}{2}, \frac{\theta}{2}}(x_0, t_0)} u\leqslant \gamma 
\bigg(\Big[\frac{1}{\theta}+\sum\limits_{i:\lambda_i=1} \frac{1}{\rho^{p_i}_i}+\sum\limits_{i:\Lambda_i=1} \frac{1}{\rho^{q_i}_i}\Big]^{\frac{N+p}{p}}\iint\limits_{Q_{\vec{\rho}, \theta}(x_0, t_0)} u^{m+1}\,dx dt\bigg)^{\frac{p}{(m+1)p+N(p|\frac{\lambda}{p}|-1)}} +\\+\gamma \bigg[\sum\limits_{i:\lambda_i<1}\Big(\frac{\theta}{\rho^{p_i}_i}\Big)^{\frac{1}{1-\lambda_i}}+
\sum\limits_{i:\Lambda_i<1}\Big(\frac{\theta}{\rho^{q_i}_i}\Big)^{\frac{1}{1-\Lambda_i}}\bigg],
\end{multline}
and finally in the degenerate case $1< \Lambda=\mathcal{L}$ we have
\begin{multline}\label{eq1.6}
\sup\limits_{Q_{\frac{\vec{\rho}}{2}, \frac{\theta}{2}}(x_0, t_0)} u\leqslant \gamma \bigg(\Big[\sum\limits_{i:\lambda_i=\Lambda} \frac{1}{\rho^{p_i}_i}+\sum\limits_{i:\Lambda_i=\Lambda} \frac{1}{\rho^{q_i}_i}\Big]^{\frac{N+p}{p}}\iint\limits_{Q_{\vec{\rho}, \theta}(x_0, t_0)} u^{m+\Lambda}\,dx dt\bigg)^{\frac{p}{(m+1)p+N(p|\frac{\lambda}{p}|-\Lambda)}} +\\+\gamma \bigg[\sum\limits_{i:\lambda_i=\Lambda}\Big(\frac{\rho^{p_i}_i}{\theta}\Big)^{\frac{1}{\Lambda-1}}+
\sum\limits_{i:\Lambda_i=\Lambda}\Big(\frac{\rho^{q_i}_i}{\theta}\Big)^{\frac{1}{\Lambda-1}}
+\sum\limits_{i:\lambda_i<\Lambda}\Big(\frac{1}{\rho^{p_i}_i}\sum\limits_{l:\lambda_l=\Lambda}\rho^{p_l}_l\Big)
^{\frac{1}{\Lambda-\lambda_i}}+\\+\sum\limits_{i:\Lambda_i<\Lambda}\Big(\frac{1}{\rho^{q_i}_i}
\sum\limits_{l:\Lambda_l=\Lambda}\rho^{q_l}_l\Big)
^{\frac{1}{\Lambda-\Lambda_i}}\bigg]\,.
\end{multline}
    
\end{remark}

\subsubsection{Case  $\mathcal{L} = \mathcal{M}$}
Let us set 
\[y_j:=\iint\limits_{Q_j}(u^{m}- k^{m}_{j})^{p_*}_{+}\,\,dx dt,\]
\[p_*=\alpha p+(1+\frac{1}{m})\frac{p}{N}=1+\frac{\mathcal{L}}{m},\quad \alpha p=1+\frac{p}{m}\Big|\frac{\lambda}{p}\Big|\,.\]
\noindent Using Lemma \ref{lem2.3} and \eqref{eq3.6}, similarly to \eqref{eq3.8} we obtain
\begin{equation}\label{eq3.10}
y_{j+1}\leqslant \gamma 2^{j\gamma}\,\big[\mathcal{H}(\theta, \vec{\rho})\big]^{\frac{N+p}{N}}\, y^{\frac{N+p}{N}}_j.
\end{equation}
By Lemma \ref{lem2.2}, inequality \eqref{eq3.10} yields $\lim\limits_{j\rightarrow \infty} y_j=0$, provided that
\begin{equation}\label{eq3.11}
y_0=\iint\limits_{Q_{\vec{\rho},\theta}}(u^{m}- \frac{k^{m}}{2})^{p_*}_{+}\,\,dx dt\leqslant \frac{1}{\gamma}\, \big[\mathcal{H}(\theta, \vec{\rho})\big]^{-\frac{N+p}{p}}.
\end{equation}
Since
$$\iint\limits_{Q_{\vec{\rho},\theta}}(u^{m}- \frac{k^{m}}{2})^{p_*}_{+}\,\,dx dt\leqslant \iint\limits_{Q_{\vec{\rho},\theta}\cap\{u\geqslant k 2^{-\frac{1}{m}}\}}u^{m p_*}\,dx dt$$
and using Chebyshev's inequality together with the absolute continuity of the integral with respect to the domain
\[\lim\limits_{k\rightarrow \infty} |[u\ge k2^{-1/m}]| \leq \lim\limits_{k\rightarrow \infty} \frac{1}{k^{mp^*}2^{-p^*}} \int_{Q_{\vec{\rho},\theta}} u^{mp^*}\, dxdt=0 \quad \Rightarrow \]
\[\quad \Rightarrow \quad \lim\limits_{k\rightarrow \infty}\iint\limits_{Q_{\vec{\rho},\theta}\cap\{u\geqslant k 2^{-\frac{1}{m}}\}}u^{m p_*}\,dx dt=0,\]
we choose $k$ large enough that \eqref{eq3.11} holds, and hence 
$$\sup\limits_{Q_{\frac{\vec{\rho}}{2}, \frac{\theta}{2}}(x_0, t_0)} u\leqslant k,$$
completing the proof of Theorem \ref{th1.1}.

\subsection{Proof of Theorem \ref{th1.2}}

\subsubsection{Case $\mathcal{L}> \mathcal{M}$.}

 Due to the fact that $\mathcal{L}> p\big|\frac{\lambda}{p}\big|+(m+1)\frac{p}{N}$ we cannot follow the proof of Theorem \ref{th1.1}. Let us estimate the  integral on the right-hand side of \eqref{eq3.6}, for this set
$$\epsilon=\frac{s-m-\mathcal{L}}{s-m p_*},\quad p_*= \alpha p+(1+\frac{1}{m})\frac{p}{N},\quad \alpha p=1+\frac{p}{m}\Big|\frac{\lambda}{p}\Big|.$$ 
Note that by assumptions $\mathcal{L}>\mathcal{M}$ and  \eqref{eq1.9} one can estimate
\begin{equation*}
s>\mathcal{L}-1+\frac{N}{p}(\mathcal{L}-p\Big|\frac{\lambda}{p}\Big|)>\mathcal{L}+m>p\Big|\frac{\lambda}{p}\Big|+m+(m+1)\frac{p}{N}=m \alpha p+(m+1)\frac{p}{N}=m p_*,
\end{equation*}
so that $\epsilon \in (0,1)$ is well defined.
Fix number $d$ such that
$$d:=\Big(\iint\limits_{Q_{\vec{\rho}, \theta}(x_0, t_0)} u^s dx dt\Big)^{\frac{1}{s}},$$ 
then by the H\"{o}lder inequality, and using the fact that $\dfrac{m+\mathcal{L}-m \epsilon p_*}{1-\epsilon}=s$, we obtain
\begin{multline}\label{eq3.12}
\iint\limits_{Q_{j}}(u^{m}- k^{m}_{j})^{1+\frac{\mathcal{L}}{m}}_{+}\,\,dx dt\leqslant \gamma \iint\limits_{Q_{j}}(u^{m}-k^{m}_{j})^{\epsilon p_*}_{+}\,u^{m+\mathcal{L}- m \epsilon p_*}\,\,dx dt \\\leqslant \Big(\iint\limits_{Q_{j}}(u^{m}-k^{m}_{j})^{p_*}_{+}\,\,dx dt\Big)^{\epsilon} \Big(\iint\limits_{Q_{j}}u^s\,\,dx dt\Big)^{1-\epsilon}\\\leqslant \gamma 2^{j\gamma} d^{s(1-\epsilon)}\,\Big(\iint\limits_{Q_j}(u^{m}-k^{m}_{j})^{p_*}_{+}\,\,dx dt\Big)^\epsilon .
\end{multline}
Inequalities \eqref{eq3.6},\eqref{eq3.12} yield
\begin{multline}\label{eq3.13}
\sup\limits_{t\in I_j}\int\limits_{B_j}(u^{m}-k^{m}_{j+1})^{1+\frac{1}{m}}_{+}\zeta^{q}_j\,dx+
\sum\limits_{i=1}^N\iint\limits_{Q_{j}}|D_i (u^{m}-k^{m}_{j+1})^{\alpha_i}_{+}|^{p_i}\zeta^{q}_j\,dx\,dt \\\leqslant
\gamma 2^{j\gamma}\,d^{s(1-\epsilon)}\,\mathcal{H}(\theta, \vec{\rho})\,\Big(\iint\limits_{Q_j}(u^{m}-k^{m}_{j})^{p_*}_{+}\,\,dx dt\Big)^\epsilon .
\end{multline}
Similarly to \eqref{eq3.10}, by Lemma \ref{lem2.3}, using the H\"{o}lder inequality and \eqref{eq3.13},  we have
\begin{equation}\label{eq3.14}
y_{j+1}:=\iint\limits_{Q_{j+1}}(u^{m}-k^{m}_{j+1})^{p_*}_{+}\,dx\,dt\leqslant \gamma 2^{j\gamma}\,\Big[d^{s(1-\epsilon)}\,\mathcal{H}(\theta, \vec{\rho})\Big]^{\frac{N+p}{N}}\,y^{\epsilon \frac{N+p}{N}}_j.
\end{equation}
Calculations give  
\begin{equation}\label{eq3.15}
\epsilon(N+p)-N=\frac{p(s-\mathcal{L}+1)+N(p\big|\frac{\lambda}{p}\big|-\mathcal{L})}{s-m p_*}=\frac{\varkappa_s}{s-m p_*}>0,
\end{equation}
and hence $\epsilon\frac{N+p}{N} >1$, so, from \eqref{eq3.14}, by Lemma \ref{lem2.2} we obtain that $\lim\limits_{j\rightarrow \infty} y_j=0$, provided that
\begin{equation}\label{eq3.16}
y_0 \leqslant \frac{1}{\gamma}\,\Big[d^{s(1-\epsilon)} \mathcal{H}(\theta, \vec{\rho})\Big]^{-\frac{N+p}{\epsilon(N+p)-N}}.
\end{equation}
By the H\"{o}lder inequality
\begin{multline*}
y_0=\iint\limits_{Q_{\vec{\rho},\theta}}(u^{m}- \frac{k^{m}}{2})^{p_*}_{+}\,\,dx dt\leqslant \iint\limits_{Q_{\vec{\rho},\theta}
\cap\{u\geqslant k 2^{-\frac{1}{m}}\}}u^{m p_*}\,dx dt\leqslant\\\leqslant  d^{mp_*}\big|Q_{\vec{\rho}, \theta}(x_0, t_0)\cap\{u\geqslant k 2^{-\frac{1}{m}}\}\big|^{1-\frac{mp_*}{s}}.
\end{multline*}
Evidently we have
$$\big|Q_{\vec{\rho}, \theta}(x_0, t_0)\cap\{u\geqslant k 2^{-\frac{1}{m}}\}\big|\leqslant \Big(\frac{2^{\frac{1}{m}}\,d}{k}\Big)^s,$$
and hence
\begin{equation}\label{eq3.17}
y_0\leqslant \frac{2^{\frac{s}{m}(1-\frac{mp_*}{s})} d^s}{k^{s-mp_*}}.
\end{equation}
Therefore, \eqref{eq3.17} implies \eqref{eq3.16}, provided that we choose $k$ large enough such that 
\begin{equation}\label{eq3.18}
k^{s-m p_*}=\gamma\, d^{s+\frac{s(1-\epsilon)(N+p)}{\epsilon(N+p)-N}}\,\big[\mathcal{H}(\theta, \vec{\rho})\big]^{\frac{N+p}{\epsilon(N+p)-N}}=\gamma d^{\frac{sp}{\epsilon(N+p)-N}}\big[\mathcal{H}(\theta, \vec{\rho})\big]^{\frac{N+p}{\epsilon(N+p)-N}}.
\end{equation}
Taking into account \eqref{eq3.15} and the definition of $d$,  \eqref{eq3.18} can be rewritten as
$$k=\gamma d^{\frac{sp}{\varkappa_s}}\,\big[\mathcal{H}(\theta, \vec{\rho})\big]^{\frac{N+p}{\varkappa_s}}=\gamma \big[\mathcal{H}(\theta, \vec{\rho})\big]^{\frac{N+p}{\varkappa_s}}\Big(\iint\limits_{Q_{\vec{\rho}, \theta}(x_0, t_0)} u^s \,dx dt\Big)^{\frac{p}{\varkappa_s}},$$
and hence
$$\sup\limits_{Q_{\frac{\vec{\rho}}{2}, \frac{\theta}{2}}(x_0, t_0)} u\leqslant k,$$
which, taking into account \eqref{eq3.4}, \eqref{eq3.5}, proves Theorem \ref{th1.2} in the case $\mathcal{L}> p\big|\frac{\lambda}{p}\big|+(m+1)\frac{p}{N}$.

\subsubsection{Case $\mathcal{L}= \mathcal{M}$.}
Set 
$$y_j:=\iint\limits_{Q_j}(u^{m}- k^{m}_{j})^{p_*}_{+}\,\,dx dt,\quad p_*=\alpha p+(1+\frac{1}{m})\frac{p}{N}=1+\frac{\mathcal{L}}{m},\quad \alpha p=1+\frac{p}{m}\Big|\frac{\lambda}{p}\Big|,$$
by \eqref{eq3.10}, \eqref{eq3.11} $\lim\limits_{j\rightarrow \infty} y_j=0$, provided that
\begin{equation}\label{eq3.19}
y_0\leqslant \frac{1}{\gamma} \big[\mathcal{H}(\theta, \vec{\rho})\big]^{-\frac{N+p}{p}},
\end{equation}
with $\mathcal{H}(\theta, \vec{\rho})$ defined in \eqref{eq3.5}. By \eqref{eq3.17}
\begin{equation*}
y_0\leqslant \frac{2^{\frac{s}{m}(1-\frac{mp_*}{s})} d^s}{k^{s-mp_*}},\quad d =\Big(\iint\limits_{Q_{\vec{\rho},\theta}(x_0, t_0)} u^s\,dx dt\Big)^{\frac{1}{s}}.
\end{equation*}
Therefore, inequality \eqref{eq3.19} holds, provided that $k$ is chosen to satisfy
$$k^{s-mp_*}=\gamma \big[\mathcal{H}(\theta, \vec{\rho})\big]^{\frac{N+p}{p}}\,\iint\limits_{Q_{\vec{\rho},\theta}(x_0, t_0)} u^s\,dx dt.$$
By our choice of $\mathcal{L}$, calculations give   $s-m\,p_*=\frac{\varkappa_s}{p}.$ Taking into account \eqref{eq3.4}, we arrive at the required \eqref{eq3.10}, \eqref{eq3.11}, this completes the proof of Theorem \ref{th1.2}.

\section{Examples}
\noindent In this Section we show that local weak sub-solutions to many evolutionary nonlinear equations are members of $\mathcal{PDG}^+(\Omega_T)$. 

\subsection{Example \eqref{eq4.1} - Doubly Nonlinear Equations.}
\label{ex4.1}
Local weak solutions to doubly nonlinear parabolic equation with standard growth \eqref{eq4.1} satisfy inequality \eqref{eq1.1} in the cylinder $Q_{r,\eta}(y, \tau)=B_r(y)\times (\tau-\eta, \tau)$ with 
\[m=m_1=...=m_N=n_1=...=n_N, \quad p=p_1=...=p_N=q_1=...=q_N, \]
see, for example \cite{BogDuzGiaLiaSch}. Conditions \eqref{eq1.3} and $\mathcal{L}>\mathcal{M}$ are translated into
$$1\leqslant m(p-1)+(m+1)\frac{p}{N},\quad \text{or} \quad 1> m(p-1)+(m+1)\frac{p}{N}.$$
The numbers $\varkappa_s$ defined in \eqref{eq1.7} can be rewritten as
$$\varkappa_s= ps+ N(m(p-1)-1) > 0.$$

\subsection{Example \eqref{eq4.2} - Generalized Orlicz Growth.} \label{ex4.2}
Doubly nonlinear parabolic equations with generalized Orlicz growth as \eqref{eq4.2}, i.e. 
\begin{equation*}
u_t- div\Big(\varphi(x, t, u, |\nabla u^{m^-}|)\frac{\nabla u^{m^-}}{|\nabla u^{m^-}|}\Big)=0,\quad (x,t) \in \Omega_T,
\end{equation*}
where $\varphi(x, t, u, v)$ is increasing in $v$ for all $(x,t)\in \Omega_T$ and $u>0$ and satisfies the conditions
\begin{equation}\label{eq4.3}
K_1 u^{(m-m^-)(p-1)} v^{p}\leqslant \varPhi(x, t, u, v):=\varphi(x, t, u, v)\,v\leqslant K_2 u^{(n-m^-)(q-1)} v^{q},
\end{equation}
where $m^-:=\min (m, n)>0$, $p\leqslant q$, $m(p-1)\leqslant n (q-1)$.

\begin{remark} The case of double-phase equations is covered by this one; the only difference being that the constants $K_1,K_2$ above and therefore the constant $C$ of the definition of $\mathcal{PDG}^+(\Omega_T)$ (see \eqref{eq1.1}) depend on the local maximum and on the minimum value of the phase (see \cite{CHS2}, Lemma 3.1).  

    \end{remark}

\begin{definition}
A function
$$\begin{cases}
u\in C_{loc}(0, T; L^{1+m^-}_{loc}(\Omega_T)),\quad u^{m^-}\in L^{q}(0, T; W^{1,q}_{loc}(\Omega)),\\
u^{\frac{m(p-1)+m^-}{p_1}}\in L^{p}_{loc}(0, T; W^{1,p}_{loc}(\Omega)),\quad u^{\frac{n(q-1)+m^-}{q}}\in  L^{q}_{loc}(0, T; L^{q}_{loc}(\Omega))
\end{cases} $$
is a local, weak subsolution to \eqref{eq4.2} if for every compact set $K\subset \Omega$ and every subinterval $[t_1, t_2]\subset (0, T]$ the inequality
\begin{equation}\label{eq4.4}
\int\limits_K u\,\zeta dx \Big|^{t_2}_{t_1}+ \int\limits^{t_2}_{t_1}\int\limits_K\big\{-u\,\zeta_t+\varphi(x, t, u, |\nabla u^{m^-}|)\frac{\nabla u^{m^-}}{|\nabla u^{m^-}|} \nabla \zeta\big\}\,dx dt\leqslant 0,
\end{equation}
is valid for all test functions $0 \leq \zeta \in W^{1, 1+m^-}_{loc}(0, T; L^{1+m^-}(K))\cap L^{q}_{loc}(0, T; W^{1, q}_0(K))$.
\end{definition} \noindent We show that local weak sub-solutions to \eqref{eq4.2} belong to $\mathcal{PDG}^+ (\Omega_T)$.
\begin{proof} First we observe that the following inequality holds true:
\begin{equation}\label{claim}
\varphi(x, t, u, a)\,b\leqslant \epsilon \varphi(x, t, u, a)\,a+ \varphi(x, t, u, \frac{b}{\epsilon})\,b,\quad (x, t)\in \Omega_T,\quad \epsilon, u, a, b >0.
\end{equation} 
Indeed, if $b\leqslant \epsilon a$, then $\varphi(x, t, u, a)\,b\leqslant \epsilon \varphi(x, t, u, a)\,a$, and if $b\geqslant \epsilon a$,
by the fact that $\varphi(\cdot, \cdot, \cdot, a)$ is increasing, we obtain $\varphi(x, t, u, a)\,b\leqslant  \varphi(x, t, u, \frac{b}{\epsilon})\,b$, from which the claim follows.
\vskip0.2cm \noindent Now, let us test \eqref{eq4.4} by $(u^{m^-}-k^{m^-})_+\,\zeta^{q}$, where $\zeta(x, t)$ is a piecewise smooth, cutoff function,
vanishing on $\partial B_r(y)\times (\tau-\eta, \tau)$, $0\leqslant \zeta(x, t)\leqslant 1$ and integrate over $B_r(y)\times (\tau-\eta, t_1)$ with $t_1 \in (\tau-\eta, \tau)$. The use of such a test function is justified, modulus a standard averaging process (see \cite{DiB} Lemma 3.2 Chap I), we obtain
\begin{multline*}
\int\limits_{B_r(y)\times\{t_1\}}g_+(u^{m^-}, k^{m^-})\,\zeta^q\,dx+
\int\limits_{\tau-\eta}^{t_1}\int\limits_{B_r(y)} \varPhi(x, t, u, |\nabla (u^{m^-}-k^{m^-})_+|)\,\zeta^q\,dx dt\\\leqslant \int\limits_{B_r(y)\times\{\tau-\eta\}}g_{+}(u^{m^-}, k^{m^-})\,\,\zeta^{q}\,dx+ \int\limits_{\tau-\eta}^{t_1}\int\limits_{B_r(y)} g_{+}(u^{m^-}, k^{m^-})\,|\zeta_t|\,dx dt+\\+
q\int\limits^{t_1}_{\tau-\eta}\int\limits_{B_r(y)}\varphi(x, t, u, |\nabla(u^{m^-}-k^{m^-})_+|)(u^{m^-}-k^{m^-})_+|\nabla \zeta|\,\zeta^{q-1} dx dt.
\end{multline*}
We estimate the last term of this inequality using \eqref{claim} with 
\[a=|\nabla(u^{m^-}-k^{m^-})_+|, \qquad b=q (u^{m^-}-k^{m^-})_+|\nabla \zeta|\,\zeta^{-1}, \quad \text{and} \quad \epsilon=\dfrac{1}{2},\]
to get
\begin{multline*}
\int\limits_{B_r(y)\times\{t_1\}}g_+(u^{m^-}, k^{m^-})\,\zeta^q\,dx+
\frac{1}{2}\int\limits_{\tau-\eta}^{t_1}\int\limits_{B_r(y)} \varPhi(x, t, u, |\nabla (u^{m^-}-k^{m^-})_+|)\,\zeta^q\,dx dt\\\leqslant \int\limits_{B_r(y)\times\{\tau-\eta\}}g_{+}(u^{m^-}, k^{m^-})\,\,\zeta^{q}\,dx+ \int\limits_{\tau-\eta}^{t_1}\int\limits_{B_r(y)} g_{+}(u^{m^-}, k^{m^-})\,|\zeta_t|\,dx dt+\\+
q\int\limits^{t_1}_{\tau-\eta}\int\limits_{B_r(y)}\varphi(x, t, u, 2 q (u^{m^-}-k^{m^-})_+||\nabla \zeta|\zeta^{-1})(u^{m^-}-k^{m^-})_+|\nabla \zeta|\,\zeta^{q-1} dx dt.
\end{multline*}
From this by \eqref{eq4.3}
\begin{multline}\label{eq4.5}
\int\limits_{B_r(y)\times\{t_1\}}g_+(u^{m^-}, k^{m^-})\,\zeta^q\,dx+
\frac{K_1}{2}\int\limits_{\tau-\eta}^{t_1}\int\limits_{B_r(y)} u^{(m-m^-)(p-1)}|\nabla(u^{m^-}-k^{m^-})_+|^p\,\zeta^q\,dx dt\\\leqslant
\int\limits_{B_r(y)\times\{\tau-\eta\}}g_{+}(u^{m^-}, k^{m^-})\,\,\zeta^{q}\,dx+ \int\limits_{\tau-\eta}^{t_1}\int\limits_{B_r(y)} g_{+}(u^{m^-}, k^{m^-})\,|\zeta_t|\,dx dt+\\+\gamma (q, K_2)\int\limits_{\tau-\eta}^{t_1}\int\limits_{B_r(y)} u^{(n-m^-)(q-1)}(u^{m^-}-k^{m^-})^q_+\,|\nabla \zeta|^q\,dx dt,
\end{multline}
that is inequality \eqref{eq1.1} with 
\[m_1= ...=m_N=m, \quad n_1= ...=n_N=n, \quad p_1= ...=p_N=p, \quad q_1= ...=q_N=q\,.\]
\noindent Conditions \eqref{eq1.3} and $\mathcal{L}> \mathcal{M}$ are translated into
$$\max(1, n(q-1))\leqslant m(p-1)+(m^-+1)\frac{p}{N},\quad \text{or} \quad \max(1, n(q-1))> m(p-1)+(m^-+1)\frac{p}{N}.$$
The numbers $\varkappa_s$ defined in \eqref{eq1.7} can be rewritten as
$$\varkappa_s= 
\begin{cases}
ps+ N(m(p-1)-1),\quad  \text{if}\quad n(q-1)<1,\\
p(s+1-n(q-1))+N(m(p-1)-n(q-1)),\quad \text{if}\quad n(q-1)>1.
\end{cases}$$

\end{proof}

\subsection{Example \eqref{eq4.6} - Doubly Nonlinear Anisotropic Equations.} \label{ex4.3}
Here we cover the case of doubly nonlinear anisotropic parabolic equations such as \eqref{eq4.6} in $\Omega_T$, with $p_i>1$ for all $i=1\dots, N$ and $m=\min(m_1, ..., m_N)>0$.
\begin{definition}
    A function
$$\begin{cases}
u\in C_{loc}(0, T; L^{1+m}_{loc}(\Omega_T)),\quad u^{m}\in L^{\vec{p}}(0, T; W^{1,\vec{p}}_{loc}(\Omega)),\\
u^{\frac{m_i(p_i-1)+m}{p_i}},  D_i \big(u^{\frac{m_i(p_i-1)+m}{p_i}}\big)\in L^{p_i}(0, T; L^{p_i}_{loc}(\Omega)),\quad i=1, ..., N,
\end{cases} $$
is a local, weak subsolution to \eqref{eq4.6} if for every compact set $K\subset \Omega$ and every subinterval $[t_1, t_2]\subset (0, T]$ the inequality
\begin{equation}\label{eq4.7}
\int\limits_K u\,\zeta dx \Big|^{t_2}_{t_1}+ \int\limits^{t_2}_{t_1}\int\limits_K\big\{-u\,\zeta_t+\sum\limits_{i=1}^N (u^{(m_i-m)(p_i-1)}|D_i u^{m}|^{p_i-2} D_i u^{m} D_i\zeta\big\}\,dx dt\leqslant 0,
\end{equation}
is valid for all non-negative testing functions $\zeta \in W^{1, 1+m}_{loc}(0, T; L^{1+m}(K))\cap L^{\vec{p}}_{loc}(0, T; W^{1, \vec{p}}_0(K))$.
\end{definition}
\noindent We show that local weak sub-solutions to \eqref{eq4.6} belong to $\mathcal{PDG}^+ (\Omega_T)$.
\begin{proof} Test \eqref{eq4.7} by $(u^{m}-k^{m})_+\,\zeta^{p^+}$, where $p^+:=\max(p_1, ..., p_N)$ and $\zeta(x, t)$ is a piecewise smooth, cutoff function, vanishing on $\partial K_{\vec{r}}(y)\times (\tau-\eta, \tau)$, $0\leqslant \zeta(x, t)\leqslant 1$ and integrate over $K_{\vec{r}}(y)\times (\tau-\eta, t_1)$ with $t_1 \in (\tau-\eta, \tau)$. The use of such a test function is justified, modulus a standard averaging process (see \cite{DiB} Lemma 3.2 Chap I), we obtain
\begin{multline*}
\int\limits_{K_{\vec{r}}(y)\times\{t_1\}}g(u^{m}, k^{m})\,\zeta^{p^+}\,dx+
\sum\limits_{i=1}^N\int\limits_{\tau-\eta}^{t_1}\int\limits_{K_{\vec{r}}(y)}u^{(m_i-m)(p_i-1)}|D_i (u^{m}-k^{m})_+|^{p_i}\,\zeta^{p^+} dx dt \\\leqslant \int\limits_{K_{\vec{r}}(y)\times\{\tau-\eta\}}g(u^{m}, k^{m})\,\,\zeta^{p^+}\,dx+ \int\limits_{\tau-\eta}^{t_1}\int\limits_{K_{\vec{r}}(y)} g(u^{m}, k^{m})\,|\zeta_t|\,dx dt+\\+
\sum\limits_{i=1}^N p_i\,\int\limits_{\tau-\eta}^{t_1}\int\limits_{K_{\vec{r}}(y)}u^{(m_i-m)(p_i-1)}|D_i (u^{m}-k^{m})_+|^{p_i-1}\,(u^{m}-k^{m})_+\,|D_i\zeta|\zeta^{p^+-1} dx dt,
\end{multline*}
from which by the Young inequality we arrive at
\begin{multline}\label{eq4.8}
\int\limits_{K_{\vec{r}}(y)\times\{t_1\}}g(u^{m}, k^{m})\,\zeta^{p^+}\,dx+
\frac{1}{2}\sum\limits_{i=1}^N\int\limits_{\tau-\eta}^{t_1}\int\limits_{K_{\vec{r}}(y)}u^{(m_i-m)(p_i-1)}|D_i (u^{m}-k^{m})_+|^{p_i}\,\zeta^{p^+} dx dt \\\leqslant \int\limits_{K_{\vec{r}}(y)\times\{\tau-\eta\}}g(u^{m}, k^{m})\,\,\zeta^{p^+}\,dx+ \int\limits_{\tau-\eta}^{t_1}\int\limits_{K_{\vec{r}}(y)} g(u^{m}, k^{m})\,|\zeta_t|\,dx dt+\\+\gamma(p_1, ..., p_N)\sum\limits_{i=1}^N \,\int\limits_{\tau-\eta}^{t_1}\int\limits_{K_{\vec{r}}(y)}u^{(m_i-m)(p_i-1)}\,(u^{m}-k^{m})^{p_i}_+\,|D_i\zeta|^{p_i}dx dt,
\end{multline}
that is inequality \eqref{eq1.1} with $m_i=n_i$, $p_i=q_i$, $i=1, ..., N$.  Conditions \eqref{eq1.3}, $\mathcal{L}>\mathcal{M}$ are translated into
$$\max(1, \Lambda)\leqslant p\Big|\frac{\lambda}{p}\Big|+(m+1)\frac{p}{N},\quad \text{or} \quad \max(1, \Lambda)> p\Big|\frac{\lambda}{p}\Big|+(m+1)\frac{p}{N},$$
$\Lambda:=\max(m_1(p_1-1), ..., m_N(p_N-1))$. The numbers $\varkappa_s$ defined in \eqref{eq1.7} can be rewritten as
$$\varkappa_s= 
\begin{cases}
ps+ N(p\big|\frac{\lambda}{p}\big|-1),\quad  \text{if}\quad \Lambda<1,\\
p(s+1-\Lambda)+N(p\big|\frac{\lambda}{p}\big|-\Lambda),\quad \text{if}\quad \Lambda>1.
\end{cases}$$
 
\end{proof}

\section{Acknowledgements}
 Igor Skrypnik is partially supported  by the grant ‘Mathematical modelling of complex systems and processes related to security’ of the National Academy of Sciences of Ukraine under the budget programme ‘Support for the development of priority areas of scientific research’ for 2025-2026, p/n 0125U000299. Simone Ciani is partially founded by GNAMPA (INdAM), and PNR italian fundings. E. Henriques was financed by Portuguese Funds through FCT - Funda\c c\~ao para a Ci\^encia e a Tecnologia - within the Projects UIDB/00013/2020 and UIDP/00013/2020 with the references DOI 10.54499/UIDB/00013/2020 (https://doi.org/ 10.54499/UIDB/00013/2020) and DOI 10.54499/UIDP/00013/2020 (https://doi.org/10.54499/ UIDP/00013/2020).

\bigskip

CONTACT INFORMATION

\medskip

{\bf Simone~Ciani}\\Università di Bologna Alma Mater, Piazza Porta San Donato 5, Italy\\
simone.ciani3@unibo.it

\medskip
{\bf Eurica Henriques}\\Centro de Matemática CMAT; Polo CMAT-UTAD
Departamento de Matemática
Universidade de Trás-os-Montes e Alto Douro, Vila Real, Portugal
  \\
eurica@utad.pt

\medskip 
{\bf Mariia Savchenko}\\ Institute of Applied Mathematics and Mechanics,
National Academy of Sciences of Ukraine, Gen. Batiouk Str. 19, 84116 Sloviansk, Ukraine\\
shan\textunderscore maria@ukr.net

\medskip

{\bf Igor I.~Skrypnik}\\Institute of Applied Mathematics and Mechanics,
National Academy of Sciences of Ukraine, Gen. Batiouk Str. 19, 84116 Sloviansk, Ukraine\\
ihor.skrypnik@gmail.com

\medskip

\medskip


\begin{thebibliography}{99}

\bibitem{BiaCupMas}
S. Biagi, G. Cupini, E. Mascolo, Regularity of quasi-minimizers for non-uniformly elliptic
integrals, J. Math. Anal. Appl., 485 (2020), 123838.

\bibitem{BocMarSbo}
L. Boccardo, P. Marcellini, C. Sbordone,  $L^\infty$ regularity for variational problems with sharp
nonstandard growth conditions, Boll. Un. Mat. Ital. A (7), 4 (1990), 219-225.

\bibitem{BonVaz}
M. Bonforte and J.L. V\'azquez, Positivity, Local Smoothing and Harnack Inequalities for Very Fast Diffusion Equations, Adv. Math., 223, (2010), 529-578.

\bibitem{BonIagVaz}
M. Bonforte, R.G. Iagar and J.L. V\'azquez, Local smoothing effects, positivity,
and Harnack inequalities for the fast $p$-Laplacian equation, Adv. Math., 224,
(2010), 2151-2215.


\bibitem{BogDuzGiaLiaSch}
V. B\"{o}gelein, F. Duzaar, U. Gianazza, N. Liao, C. Scheven,
H\"{o}lder Continuity of the Gradient of Solutions to Doubly Non-Linear Parabolic Equations, Preprint, https://doi.org/10.48550/arXiv.2305.08539.


 
\bibitem{CarGaoGioLeo}
M. Carozza, H. Gao, R. Giova, F. Leonetti, A boundedness result for minimizers
of some polyconvex integrals, J. Optim. Theory Appl., 178 (2018), 699–725.

\bibitem{Cia}
A. Cianchi, Local boundedness of minimizers of anisotropic functionals, Ann. Inst. H. Poincar\'e C
Anal. Non Lin\'eaire, 17 (2000), 147–168.


\bibitem{CH} S. Ciani, E. Henriques, Harnack-type estimates and extinction in finite time for a class of anisotropic porous medium type equations. Calculus of Variations and Partial Differential Equations, 64(4), (2025), 117.

\bibitem{CHS} S. Ciani, E. Henriques, I.I. Skrypnik, The impact of intrinsic scaling on the rate of extinction for anisotropic non-Newtonian fast diffusion, Nonlinear Analysis, 242, (2024), 113497.

\bibitem{CHS2} S. Ciani, E. Henriques,I.I. Skrypnik, Fine Boundary Continuity for Degenerate Double-Phase Diffusion, Potential Anal (2025), Online First https://doi.org/10.1007/s11118-025-10198-0.

\bibitem{CupMarMas}
G. Cupini, P. Marcellini, E. Mascolo, Regularity of minimizers under limit growth conditions, Nonl. Analysis, 153 (2017), 294-310.

\bibitem{CupLeoMas}
G. Cupini, F. Leonetti, E. Mascolo, Local boundedness for minimizers of some polyconvex
integrals, Arch. Rational Mech. Anal., 224 (2017), 269–289.

\bibitem{DeG}
E. De Giorgi, Sulla differenziabilit\'a e l'analiticit\'a delle estremali degli integrali multipli regolari, Mem. Accad. Sci. Torino
Cl. Sci. Fis. Mat. Natur., 3 (1957), 25-43.

\bibitem{DiB-Annali} E. Di Benedetto, On the local behaviour of solutions of degenerate parabolic equations with measurable coefficients, Annali della Scuola Normale Superiore di Pisa-Classe di Scienze, 13(3), (1986), 487-535.

\bibitem{DiB}
E. DiBenedetto, Degenerate Parabolic Equations, Universitext, Springer-Verlag, New York, 1993.

\bibitem{DiBGiaVes}
E. DiBenedetto, U. Gianazza, V. Vespri, Harnack's Inequality for Degenerate and
Singular Parabolic Equations, Springer Monographs in Mathematics, 2012.

\bibitem{Porzio2} G. Di Blasio, M.M. Porzio, Uniqueness, regularity and behavior in time of the solutions to nonlinear anisotropic parabolic equations, J. Evol. Equ. 25, (2025), 45. 

\bibitem{DMV} F.G. D\"uzg\"un, S.J.N. Mosconi, V. Vespri, Anisotropic Sobolev embeddings and the speed of propagation for parabolic equations. Journal of Evolution Equations, 19, (2019), 845-882.

 
\bibitem{FusSbo1}
N. Fusco, C. Sbordone, Local boundedness of minimizers in a limit case, Manuscripta Math., 69
(1990), 19–25.
 
\bibitem{FusSbo2}
N. Fusco, C. Sbordone, Some remarks on the regularity of minima of anisotropic integrals,
Commun. Part. Diff. Eq., 18 (1993), 153–167.

\bibitem{Gia}
M. Giaquinta, Growth conditions and regularity, a counterexample, Manuscripta Math., 59 (1987), 245-248. 

\bibitem{Giu}
E. Giusti, Direct Methods in the Calculus of Variations, World Scientific Publishing Co., Inc., River Edge, NJ, 2003.


\bibitem{Henriques1} E. Henriques, R. Laleoglu, Boundedness for Some Doubly Nonlinear Parabolic Equations in Measure Spaces, Journal of Dynamics and Differential Equations, 30, (2018), 1029-1051.

\bibitem{Henriques2} E. Henriques, R. Laleoglu, Local and global boundedness for some nonlinear parabolic equations exhibiting a time singularity, Differential and Integral Equations,
29, 11/12, (2016), 1029-1048.


\bibitem{Iva}
A. V. Ivanov, Maximum modulus estimates for generalized solutions to doubly nonlinear parabolic equations,
English translation in J. Math. Sci. (New York), 87 (1997), no. 2, 3322-3342.




\bibitem{Kol1}
I. M. Kolodii, The boundedness of generalized solutions of elliptic differential equations,
Moskow University Vestnik, 25 (1970), 44–52 (In Russian).

\bibitem{Kol2}
I. M. Kolodii,  On boundedness of generalized solution of parabolic differential equations. Moscow University Vestnik, 5,  (1971), 25-31 (In Russian).

\bibitem{LadUra}
O.A. Ladyzenskaya and N.N. Ural'tzeva, Linear and Quasilinear Elliptic Equations, Academic Press, New York, 1968.

\bibitem{LadSolUra}
O.A. Ladyzenskaya, N.A. Solonnikov and N.N. Ural'tzeva, Linear and Quasilinear Equations of Parabolic Type, Translations of Mathematical Monographs, 23, American Mathematical Society, Providence, RI, 1967.

\bibitem{Mar}
P. Marcellini, Un example de solution discontinue d'un probleme variationnel dans le cas scalaire. Preprint 11, Istituto
Matematico "U. Dini", Universit\'a di Firenze, 1987.

\bibitem{MinXit}
Y. Mingqi, L. Xiting, Boundedness of solutions of parabolic equations with anisotropic growth conditions. Canadian
Journal of Mathematics, 49(4), (1997), 798-809.


\bibitem{Ok}
J. Ok,
Regularity for double phase problems under additional integrability assumptions,
Nonlinear Anal. \textbf{194} (2020) 111408.



\bibitem{Porzio} M.M. Porzio, Regularity and time behavior of the solutions to weak monotone parabolic equations, J. Evol. Equ. 21, (2021), 3849–3889. 
\bibitem{Sin}
T. Singer, Local boundedness of variational solutions to evolutionary problems with nonstandard growth. NoDEA,
23 (2), (2016), 1-23.

\bibitem{Str}
B. Stroffolini, Global boundedness of solutions of anisotropic variational problems, Boll. Un. Mat.
Ital. A (7), 5 (1991), 345–352.

\bibitem{Tro}
M. Troisi, Teoremi di inclusione per spazi di Sobolev non isotropi, Ricerche Mat. 18 (1969) 3-24.

\bibitem{VV} V. Vespri, M. Vestberg, An extensive study of the regularity of solutions to doubly singular equations, Advances in Calculus of Variations, 15(3), (2022), 435-473.

\bibitem{Wieser} W. Wieser, Parabolic Q-minima and minimal solutions to variational flow. Manuscripta Mathematica, 59(1), (1987), 63-107.
\end{thebibliography}
\end{document}